\documentclass[reqno,10pt]{amsart}	
\def\Title{On the Statistics of Lattice Polytopes}
\def\Abstract{%
	We use the notions of reflexivity and of reflexive dimensions in order
	to introduce probability measures for lattice polytopes and 
	initiate the investigation of their statistical properties.
	Examples of applications to discrete geometry include a study of
	randomness of self-duality of reflexive polytopes and implications 
	for expectation values of the numbers of such polytopes in higher 
	dimensions. 
	We also discuss enumeration problems and related 
	algorithms and point out interesting open problems. In this context 
	we define the notion of IP-confined polytopes. Our new results 
	include the list of IP-simplices in 3 dimensions that are not 
	IP-confined.  
	The main motivation for the study of these issues co\-mes from 
	applications in algebraic geometry and string theory.
}
\def\Author{Maximilian Kreuzer}
\def\Address{Institute for Theoretical Physics, Vienna University of 
	Technology\\	Wiedner Hauptstrasse 8-10, 1040 Vienna, Austria}
\def\Email{maximilian.kreuzer@tuwien.ac.at}

\usepackage{amsmath,array,float,amsfonts,amssymb,latexsym,eufrak,cite}	
\def\url#1{#1}	

\newtheorem{theorem}{Theorem}[section]
\newtheorem{corollary}[theorem]{Corollary}

\newtheorem{code}[theorem]{Program code}	     \theoremstyle{definition}

\newtheorem{definition}[theorem]{Definition}

\newtheorem{remark}[theorem]{Remark}

\newtheorem{question}[theorem]{Question}
	    \numberwithin{equation}{section}

\usepackage{epic,pst-plot,alltt}		 \usepackage[latin1]{inputenc}

\usepackage{midfloat}
\def\ifundefined#1{\expandafter\ifx\csname#1\endcsname\relax}

\def\putlab#1)#2#3{\put#1){\makebox(0,0)[#2]{\small #3}}}

\def\IN{{\mathbb N}}	\def\IR{{\mathbb R}}	\def\IZ{{\mathbb Z}}

\def\IP{\hbox{I\hspace{-.7pt}P}}		\long\def\del#1\enddel{}

         \let\th=\theta  
           
            \let\p=\pi      	 
               
               \let\S=\Sigma 
          \let\G=\Gamma   \let\D=\Delta
					
\let\0=\over	\let\6=\partial	\def\2{{\frac12}}  \long\def\keep#1\endkeep{}
					  
\def\BE {\begin{equation}}      \def\EE {\end{equation}}        
\def\BEA{\begin{eqnarray}}      \def\EEA{\end{eqnarray}} 
\def\BP{\begin{picture}} 	\def\EP{\end{picture}}
\def\BI{\begin{itemize}} 	\def\EI{\end{itemize}}
\def\BC{\begin{center}} 	\def\EC{\end{center}}

     \def\HS#1 {\hspace*{#1pt}}      
\def\mao#1{\mathop{\rm #1}\nolimits}		\def\CH{\mao{ConvHull}}

\def\fig #1.{{fig.}\,\ref{#1}}		\def\tab #1.{{table}\,\ref{#1}}

\begin{document}   \columnsep=16pt   \pagestyle{empty}   \thispagestyle{empty}	
\ifundefined{arXiv}
	\author[M. Kreuzer]{\Author}	\address{\Address}	\email{\Email}
	\date{}	  

\title[\Title]{\Title}	\begin{abstract}\Abstract\end{abstract}
\twocolumn[\maketitle]				

\else \author[]{Maximilian Kreuzer}\address{}\email{}
	\title[]{\Title}\pagestyle{empty}

\newfont{\XLbf}{cmbx10 scaled 2800}   	\newfont{\XL}{cmr10 scaled 2800}
\newfont{\Ftit}{cmbx10 scaled 2000}	\newfont{\Fsec}{cmbx10 scaled 1200}
\newfont{\Faut}{cmbx10 scaled 1100}	\newfont{\Fadd}{cmr10 scaled 1100}
\def\Fabs{\fontfamily{phv}\selectfont}	

\twocolumn[\BP(0,0)(0,10)\psset{unit=1truemm}
			\psline(0,-14)(170,-14)(170,-34)
			\psline(0,-14)(0,-34)(170,-34)
\EP\BC
			{\Fabs\huge\bf\Ftit	\Title
\\[10mm]
			\large\bf\Faut	\Author
\\[1mm]			\large\rm\Fadd	\Address}
\bigskip\bigskip\bigskip\bigskip
\EC]

\noindent{\Fabs Abstract -- \Abstract}

\fi

\def\Section#1{\noindent{\addtocounter{section}1\Fsec \thesection. #1}}

							\vspace{-6mm}
\vspace{35pt}

\def\TransTM(#1,#2){	\psset{linecolor=red,linewidth=.9pt}
			\put(0,0){\pscircle*(\Pold,\Mold){60} }
	\Pnew=#1 \advance\Pnew by -\LogExp	\divide \Pnew by \Xfac	
	\Mnew=#2 \advance\Mnew by -\LogExp	\divide \Mnew by \Yfac
			\put(0,0){\psline(\Pold,\Mold)(\Pnew,\Mnew)}
	\Mold=\Mnew	\put(0,0){\pscircle*(\Pnew,\Mold){60} }
}	
\def\TransSS(#1,#2){	\psset{linecolor=blue,linewidth=.9pt}
			\put(0,0){\pscircle(\Pold,\Sold){60} }
	\Pnew=#1 \advance\Pnew by -\LogExp	\divide \Pnew by \Xfac	
	\Snew=#2 \advance\Snew by -\LogExp	\divide \Snew by \Yfac
			\put(0,0){\psline(\Pold,\Sold)(\Pnew,\Snew)}
	\Pold=\Pnew \Sold=\Snew	\put(0,0){\pscircle*(\Pold,\Sold){60} }
}

\long\def\TMstatistic{

\newcount\Xfac	\newcount\Yfac			\vspace*{1mm}

\unitlength=1.05mm 	\psset{unit=0.0105mm} 	\Xfac=1		\Yfac=2
			
\newcount\LogExp \LogExp=29934 
\newcount\XA \newcount\XB \newcount\XC \newcount\YA \newcount\YB \newcount\YC 
\XA=16118 	\XB=23026	\XC=\LogExp	
\YA=27631 	\YB=\LogExp	\YC=32236	

\advance\XA by -\LogExp \advance\XB by -\LogExp	\advance\XC by -\LogExp
\divide\XA by \Xfac	\divide\XB by \Xfac	\divide\XC by \Xfac

\advance\YA by -\LogExp \advance\YB by -\LogExp	\advance\YC by -\LogExp
\divide\YA by \Yfac	\divide\YB by \Yfac	\divide\YC by \Yfac

\newcount\YD	\YD=34539     \advance\YD by -\LogExp	   \divide\YD by \Yfac
\newcount\XlabU	\XlabU=2180	\def\XlabO{2450}

\def\RawDataTM{\parbox{15cm}{\footnotesize\begin{tabular}{|l|l|l|}\hline
zztm36: 12447467  74m+4s   230544737b &pp/2m=1.01934e+12 &pp/ss=9.68371e+12\\
zztm37: 17214077  74m+4s   324245690b &pp/2m=1.9495e+12  &pp/ss=1.85203e+13\\
zztm38: 26375307  179m+15s 504471363b &pp/2m=1.86503e+12 &pp/ss=3.09181e+12\\
zztm39: 32362083  189m+21s 625699081b &pp/2m=2.62483e+12 &pp/ss=2.37484e+12\\
zztm40: 41416687  224m+26s 813136082b &pp/2m=3.61887e+12 &pp/ss=2.53749e+12\\
zztm41: 54194432  224m+26s 1094626401b &pp/2m=6.19629e+12 &pp/ss=4.34473e+12\\
zztm42: 66834339  226m+27s 1359871259b &pp/2m=9.32533e+12 &pp/ss=6.12734e+12\\
zztm43: 82494748  343m+38s 1696685324b &pp/2m=9.39971e+12 &pp/ss=4.71287e+12\\
zztm44: 123132916  474m+43s 2546737882b &pp/2m=1.52994e+13 &pp/ss=8.19995e+12\\
zztm45: 141618155  545m+43s 2983011060b &pp/2m=1.77014e+13 &pp/ss=1.08468e+13\\
zztm46: 188755429  891m+48s 4029018095b &pp/2m=1.94692e+13 &pp/ss=1.54638e+13\\
zztm47: 232800688  1461m+70s 5018455639b &pp/2m=1.81137e+13 &pp/ss=1.10604e+13\\
zztm48: 286656341  1639m+70s 6355775337b &pp/2m=2.45436e+13 &pp/ss=1.67698e+13\\
zztm49: 423407680  1824m+73s 9449791165b &pp/2m=4.8179e+13 &pp/ss=3.36412e+13\\
zztm50: 490453816  3745m+112s 11010401691b &pp/2m=3.16423e+13 &pp/ss=1.91761e+13\\
zztm51: 876006696  6643m+153s 19677908185b &pp/2m=5.71016e+13 &pp/ss=3.27817e+13\\
zztm52: 959816636  7140m+161s 21800580062b &pp/2m=6.37939e+13 &pp/ss=3.55406e+13\\
zztm53: 1212578468  8761m+172s 28025544471b &pp/2m=8.30986e+13 &pp/ss=4.97007e+13\\
zztm54: 1417466160  18936m+225s 33163368399b &pp/2m=5.27393e+13 &pp/ss=3.96881e+13\\
zztm55: 1581357343  24621m+236s 37382543880b &pp/2m=5.05415e+13 &pp/ss=4.48989e+13\\
zztm56: 1842170300  25756m+238s 43460009502b &pp/2m=6.55767e+13 &pp/ss=5.99109e+13\\
zztm57: 2374020049  27249m+252s 56972851802b &pp/2m=1.0294e+14 &pp/ss=8.87499e+13\\
zztm58: 2959231747  32299m+253s 72550837548b &pp/2m=1.35033e+14 &pp/ss=1.3681e+14\\
zztm59: 3490120614  34833m+261s 86484670445b &pp/2m=1.74195e+14 &pp/ss=1.78813e+14\\
zztm60: 3918829821  35665m+262s 98842130419b &pp/2m=2.1451e+14 &pp/ss=2.23723e+14\\
zztm61: 4695402670  44170m+278s 121053492946b &pp/2m=2.48785e+14 &pp/ss=2.8527e+14\\
zztm62: 5296311764  56651m+296s 137738169193b &pp/2m=2.46931e+14 &pp/ss=3.20157e+14\\
zztm63: 6176698162  62318m+308s 162667394027b &pp/2m=3.0535e+14 &pp/ss=4.02172e+14\\
zztm64: 8145615178  79261m+351s 218136622925b &pp/2m=4.17636e+14 &pp/ss=5.38559e+14\\
zztm65: 9025107076  86323m+354s 238861422022b &pp/2m=4.70824e+14 &pp/ss=6.49977e+14\\
\hline\end{tabular}}}

\BC\BP(99,36)(-114,-10)		
	\newcount\Pold	\newcount\Mold	\newcount\Sold	
	\newcount\Pnew	\newcount\Mnew	\newcount\Snew	
	\Pold=16337	\Mold=27650	\Sold=29901	
	\advance\Pold by -\LogExp 	\advance\Mold by -\LogExp 
					\advance\Sold by -\LogExp
	\divide\Pold by \Xfac	\divide\Mold by \Yfac \divide\Sold by \Yfac


\TransTM(16337,27650)		\TransSS(16337,29901)		
\TransTM(16661,28299)		\TransSS(16661,30550)		
\TransTM(17088,28254)		\TransSS(17088,28760)		
\TransTM(17292,28596)		\TransSS(17292,28496)		
\TransTM(17539,28917)		\TransSS(17539,28562)		
\TransTM(17808,29455)		\TransSS(17808,29100)		
\TransTM(18018,29864)		\TransSS(18018,29444)		
\TransTM(18228,29872)		\TransSS(18228,29181)		
\TransTM(18629,30359)		\TransSS(18629,29735)		
\TransTM(18769,30505)		\TransSS(18769,30015)		
\TransTM(19056,30600)		\TransSS(19056,30370)		
\TransTM(19266,30528)		\TransSS(19266,30034)		
\TransTM(19474,30831)		\TransSS(19474,30451)		
\TransTM(19864,31506)		\TransSS(19864,31147)		
\TransTM(20011,31086)		\TransSS(20011,30585)		
\TransTM(20591,31676)		\TransSS(20591,31121)		
\TransTM(20682,31787)		\TransSS(20682,31202)		
\TransTM(20916,32051)		\TransSS(20916,31537)		
\TransTM(21072,31596)		\TransSS(21072,31312)		
\TransTM(21182,31554)		\TransSS(21182,31435)		
\TransTM(21334,31814)		\TransSS(21334,31724)		
\TransTM(21588,32265)		\TransSS(21588,32117)		
\TransTM(21808,32537)		\TransSS(21808,32550)		
\TransTM(21973,32791)		\TransSS(21973,32817)		
\TransTM(22089,32999)		\TransSS(22089,33041)		
\TransTM(22270,33148)		\TransSS(22270,33284)		
\TransTM(22390,33140)		\TransSS(22390,33400)		
\TransTM(22544,33352)		\TransSS(22544,33628)		
\TransTM(22821,33666)		\TransSS(22821,33920)		
\TransTM(22923,33786)		\TransSS(22923,34108)		
\psset{linecolor=black,linewidth=0.9pt}
\psline(\XA,\YA)(\XC,\YA)(\XC,\YD)(\XA,\YD)(\XA,\YA)
\psline(\XA,\YB)(\XC,\YB)\psline(\XA,\YC)(\XC,\YC)
\psset{linewidth=1.6pt}
\psline(-200,0)(200,0)		\putlab(4,0)l{$10^{13}$}
\psline(-200,\YC)(200,\YC)	\putlab(4,12)l{$10^{14}$}
\psline(-200,\YD)(200,\YD)	\putlab(4,24)l{$10^{15}$}

\psline[linewidth=1.5pt]{->}(\XA,\YA)(600,\YA) \putlab(5,-13.5)l{$p$}
\psline[linewidth=1.5pt]{->}(\XA,\YA)(\XA,2500) 
\putlab(-125,15)r{\red$p^2/2m$} \putlab(-126.5,20)r{\blue$p^2/s^2$}
\black					
\psline[linewidth=1.5pt](0,\XlabU)(0,\XlabO)\putlab(0,26)b{$10^{13}$}
\psline[linewidth=1.5pt](-7000,\XlabU)(-7000,\XlabO)\putlab(-70,26)b{$10^{10}$}
\psline[linewidth=1.5pt](-2333,\XlabU)(-2333,\XlabO)
\psline[linewidth=1.5pt](-4667,\XlabU)(-4667,\XlabO)
\psline[linewidth=1.5pt](-9333,\XlabU)(-9333,\XlabO)
\psline[linewidth=1.5pt](-11667,\XlabU)(-11667,\XlabO)
\EP\EC
}


\newcount\TDYfac	\newcount\TDZfac	
\newcount\XYfac	\newcount\XZfac	\XYfac=10	\XZfac=8
\newcount\Xaux	\newcount\Yaux	\newcount\Zaux	\newcount\YA	\newcount\ZA
\def\makeYZ(#1,#2,#3){
	\Xaux=#1 \multiply \Xaux by \XYfac \Yaux=#2 \multiply \Yaux by \TDYfac 
	\advance\Yaux by -\Xaux
	\Xaux=#1 \multiply \Xaux by \XZfac \Zaux=#3 \multiply \Zaux by \TDZfac 
	\advance\Zaux by -\Xaux
}
\def\xyzline(#1,#2,#3)(#4,#5,#6){\makeYZ(#1,#2,#3) \YA=\Yaux \ZA=\Zaux 
	\makeYZ(#4,#5,#6)	\drawline(\YA,\ZA)(\Yaux,\Zaux)
}	

\def\TwoDimRPdefs{  \linethickness{2pt} \def\TwoDimPtGrid##1##2{{\matrixput
			(-4,0)(\TDYfac,0){##1}(0,\TDZfac){##2}{\circle*\CSi}}}
\def\Aa{\BP(60,60)(-10,-8)\blue
	\xyzline(0,0,0)(0,0,3)\xyzline(0,0,0)(0,3,0)\xyzline(0,0,3)(0,3,0)
	\black	\HS-1 	\TwoDimPtGrid44 \EP}
\def\Ab{\BP(60,60)(-10,-8)\black
	\xyzline(0,0,0)(0,1,2)\xyzline(0,0,0)(0,2,1)\xyzline(0,1,2)(0,2,1)
	\black	\HS-1	\TwoDimPtGrid44	\EP}
\def\Ba{\BP(100,60)(-10,-5)\blue	
	\xyzline(0,0,0)(0,4,0)\xyzline(0,0,0)(0,2,2)\xyzline(0,2,2)(0,4,0)
	\black	\TwoDimPtGrid53	\EP}
\def\Bb{\BP(60,60)(-10,-5) \black
 	\xyzline(0,0,0)(0,2,0)\xyzline(0,0,0)(0,1,2)\xyzline(0,2,0)(0,1,2)
	\black	\TwoDimPtGrid33	\EP}
\def\Ca{\BP(60,60)(-10,-5) \blue		    \xyzline(0,0,0)(0,0,2)
	\xyzline(0,0,0)(0,2,0)\xyzline(0,2,0)(0,2,2)\xyzline(0,0,2)(0,2,2)
	\black	\TwoDimPtGrid33	\EP}
\def\Cb{\BP(60,60)(-10,-5) \black		    \xyzline(0,0,1)(0,1,0)
	\xyzline(0,2,1)(0,1,2)\xyzline(0,0,1)(0,1,2)\xyzline(0,1,0)(0,2,1)
	\black	\TwoDimPtGrid33	\EP}
\def\Da{\BP(60,60)(-10,-5) \black \xyzline(0,0,0)(0,0,2)\xyzline(0,0,0)(0,2,0)
	\xyzline(0,2,0)(0,2,1)\xyzline(0,0,2)(0,1,2)\xyzline(0,1,2)(0,2,1)
	\black	\TwoDimPtGrid33	\EP}
\def\Db{\BP(60,60)(-10,-5) \black \xyzline(0,0,1)(0,0,0)\xyzline(0,0,0)(0,1,0)
	\xyzline(0,2,1)(0,1,2)\xyzline(0,0,1)(0,1,2)\xyzline(0,1,0)(0,2,1)
	\black	\TwoDimPtGrid33	\EP}
\def\Ea{\BP(80,60)(-10,-5) 			    \xyzline(0,0,0)(0,0,2)
	\xyzline(0,0,0)(0,3,0)\xyzline(0,0,2)(0,1,2)\xyzline(0,1,2)(0,3,0)
	\black	\TwoDimPtGrid43	\EP}
\def\Eb{\BP(60,60)(-10,-5)
 	\xyzline(0,0,0)(0,1,0)\xyzline(0,2,1)(0,1,2)
	\xyzline(0,0,0)(0,1,2)\xyzline(0,1,0)(0,2,1)
	\black	\TwoDimPtGrid33	\EP}
\def\Fa{\BP(80,60)(-10,-5)			    \xyzline(0,0,0)(0,0,1)
	\xyzline(0,0,0)(0,3,0)\xyzline(0,3,0)(0,1,2)\xyzline(0,0,1)(0,1,2)
	\black	\TwoDimPtGrid43	\EP}
\def\Fb{\BP(60,60)(-10,-5)
 	\xyzline(0,1,2)(0,0,0)\xyzline(0,0,0)(0,2,0)\xyzline(0,2,0)(0,2,1)
	\xyzline(0,2,1)(0,1,2)
	\black	\TwoDimPtGrid33	\EP}

\def\SA{\BP(80,60)(-10,-5)
\red 	\xyzline(0,0,0)(0,1,2)\xyzline(0,0,0)(0,3,0)\xyzline(0,1,2)(0,3,0)
	\black	\TwoDimPtGrid43	\EP}

\def\SB{\BP(60,60)(-10,-5) 
\red 	\xyzline(0,0,0)(0,0,2)\xyzline(0,0,0)(0,2,0)
	\xyzline(0,2,0)(0,1,2)\xyzline(0,0,2)(0,1,2)
	\black	\TwoDimPtGrid33	\EP}

\def\SC{\BP(60,60)(-10,-5) 
\red 	\xyzline(0,0,0)(0,0,1)\xyzline(0,0,0)(0,2,0)
	\xyzline(0,2,0)(0,2,1)\xyzline(0,0,1)(0,1,2)\xyzline(0,1,2)(0,2,1)
	\black	\TwoDimPtGrid33	\EP}

\def\SD{\BP(60,60)(-10,-5) 
\red 	\xyzline(0,1,0)(0,0,1)\xyzline(0,1,0)(0,2,0)\xyzline(0,0,2)(0,0,1)
	\xyzline(0,0,2)(0,1,2)\xyzline(0,2,0)(0,2,1)\xyzline(0,1,2)(0,2,1)
	\black	\TwoDimPtGrid33	\EP}
}

\def\TwoDimRP{{	\begin{center}	\thicklines	\vspace{-18pt}	\TDZfac=16
		\TDYfac=\TDZfac			\TwoDimRPdefs
\unitlength=.66pt \def\CSi{6}	\Aa \HS20 \HS-7 \Ab \HS10 \HS5
\unitlength=1pt  \def\CSi{4}	\Ba\HS-15 \Bb \HS5 \Ca\HS-9 \Cb \\[-12pt]
\HS-3	\Da\HS-7 \Db\HS1  \Ea\HS-12  \Eb\HS5  \Fa\HS-12 \Fb	\\[-12pt]
	\SA \HS5   \SB \HS5 	\SC \HS5 \SD	\vspace{-7pt}	\end{center}}}


\Section{\noindent\Fsec Introduction\label{INTRO}}

Let us consider $d$-dimensional convex polytopes $\D$, which are convex 
hulls of a finite numbers of points in $\IR^d$.
Obvious notions of randomness in terms of random selections 
of coordinates of the vertices prefer simplicial polytopes, whose 
duals are simple \cite{Z}. Probability measures based on such concepts
hence are not invariant under the elementary duality exchanging vertices 
with bounding hyperplanes.

If we focus on lattice polytopes, whose vertices belong to a fixed
lattice $\IZ^d\subset\IR^d$, we are usually interested in affine unimodular 
equivalence classes, whose representatives differ by
lattice automorphisms. We will need some special classes of such polytopes:
\begin{definition}
	An \IP-polytope is a lattice polytope $\D$ with exactly one
	interior lattice point, 
	which we can choose to be the origin $0\in\IZ^d$. 
\\
	The polar (or dual) polytope $\D^*$ of a polytope $\D\subset \IR^d$ 
 	with 0 in its interior is defined by
\BE						\label{polar}
	\D^*=\{ y\in\IZ^d~\;
			|~\;\vec y\cdot\vec x\ge-1\quad\forall x\in\D\}.
\EE
	A reflexive polytope is an \IP-polytope $\D$ for which $\D^*$ is 
	also a lattice polytope.
\end{definition}

\begin{remark}The number of reflexive polytopes with given dimension is 
	finite. Because eq.\,(\ref{polar}) defines an involution, the polar 
	polytope $\D^*$ is reflexive if and only if $\D$ is reflexive.
	Reflexivity of an IP polytope is equivalent to the property that
	all facets have lattice distance 1 from the origin, i.e. there are
	no parallel lattice hyperplanes between the origin and the supporting 
	hyperplanes of the facets.
\end{remark}
\begin{definition} Haase and Melnikov \cite{HM06} have shown that every
	lattice polytope $\D$ is (isomorphic to) a face of a reflexive 
	polytope. The reflexive dimensions $\mao{rd}(\D)$ is the smallest 
	dimension for which such a reflexive polytope exists.
\end{definition}
\vspace{-3pt}
Since the number of reflexive polytopes is finite we can define probability
measures on the set of all polytopes that are uniform in fixed dimension and 
thus respect  
the polar duality of reflexive polytopes exchanging vertices
and facets. Similarly, we can use the fact that the number of polytopes $\D$
with fixed reflexive dimension $\mao{rd}(\D)$ is finite to define measures
on the set of all lattice polytopes. 
The numbers of polytopes with $\mao{rd}(\D)\le4$ are given in 
\tab tab:RD. below.

The content of this paper is organized as follows. In section 2 we summarize
results concerning the enumeration of reflexive polytopes and discuss some
relevant algorithms. In section 3 we extend the polarity of reflexive 
polytopes to an involutive duality of 
\IP-confined polytopes 
in $d>4$ dimensions and discuss 
various enumeration problems.
Section 4 contains statistical considerations and observations,
which are used in an attempt at estimating the number of reflexive 
polytopes in 5 dimensions.

{\it Acknowledgements.} The author would like to acknowledge discussions
with Benjamin Nill and support by the Austrian Research Fund FWF 
under grant no. P18679-N16.

\bigskip
\Section{Reflexivity 
	and weight vectors}	

The classification of reflexive polytopes in two dimensions can be 
done by hand on a piece of paper with the result shown in \fig fig:2d..
The interest in explicit enumerations in higher dimensions is mainly
motivated by the relation between lattice polytopes and algebraic
geometry, and in particular by the result of Batyrev \cite{B} that 
generic hypersurfaces in toric varieties are Calabi--Yau if and only if
the toric ambient space is defined in terms of a reflexive polytope.
On the algebro-geometric side of this correspondence the polar duality 
corresponds to an exchange of complex structure and K\"ahler moduli with
important applications to enumerative geometry \cite{Can,CK}. In terms 
of string theory, this amounts to determining instanton corrections from
a duality called mirror symmetry, which exchanges matter with anti-matter 
in particle physics \cite{Can}.

\begin{figure*}[htbp]
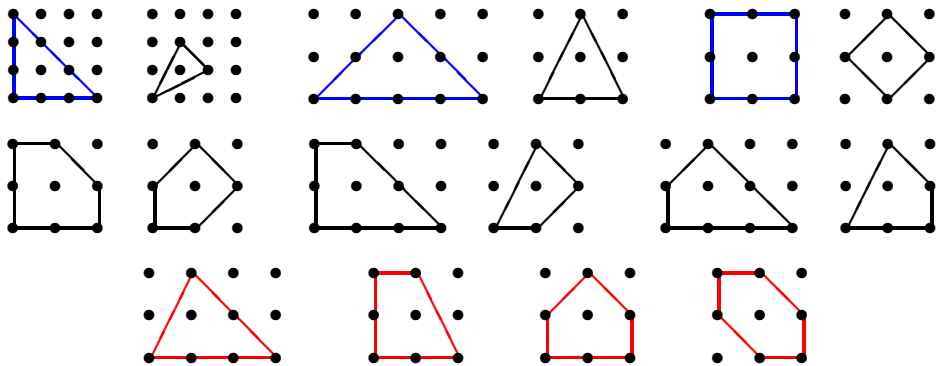
		
				\TwoDimRP
\caption{The 16 reflexive polygons in 2 dimensions:
	3 are maximal and the last 4 are selfdual.\label{fig:2d}\HS-249 }
\end{figure*}

The algorithm that has been used for the complete enumeration of reflexive 
polytopes in 3 and 4 dimensions \cite{c3d,c4d} 
is based on the idea to determine the maximal objects in terms of the polar 
minimal ones \cite{crp} and on an algorithm for the enumeration of all 
possible linear relations among vertices of the minimal polytopes \cite{S}. 
For the case of two
dimensions the three maximal reflexive polytopes are drawn, together with
the dual minimal ones, in the top row of \fig fig:2d. and it is easy to 
check that all others are subpolytopes thereof. 
If the maximal polytope is simplicial then the coefficients in the linear 
relation between the vertices of the dual minimal lattice polytope, 
which are chosen positive integral without common divisor, are called 
weight system or weight vector $\vec w=(w_0,\ldots,w_d)$. For the 
two simplicial minimal polytopes in \fig fig:2d. the weight vectors are
$(1,1,1)$ and $(1,1,2)$ and they can be obtained as $w_i=D\,q_i$ where $D$
is the common denominator of the barycentric coordinates $q_i$ of the origin
with respect to the vertices of the minimal simplices, 
i.e. $\sum q_jv_j=0$ and $\sum q_j=1$. The polar maximal polytopes can be
recovered as Newton polytopes (i.e. exponent vectors of the monomials) 
of generic quasi-homogeneous polynomials of degree $D$ w.r.t. the 
weights $w_i$. The 10 points of the
first maximal polytope, for example, correspond to the 10 monomials of
degree 3 in 3 variables.

Non-simplicial minimal \IP-polytopes are convex hulls of a 
collection of lower-dimensional minimal \IP-simplices
and thus define combined weight systems, or weight matrices, whose rows 
consist of weight vectors augmented by zeros for the vertices that do not 
belong to the respective simplex. For the third minimal polytope in 
\fig fig:2d., i.e. the square, the matrix of linear relations reads 
$1~0~1~0\choose0~1~0~1$, while a minimal pyramid in three dimensions is 
described by a weight matrix of the form
$	w_{11}~w_{12}~w_{13}\HS8 0\HS14 0\HS4	\choose\! 
		w_{21}\HS8 0\HS14	0\HS8 w_{24}~w_{25}\!$.

\begin{definition} For $w_j\in\IN$ with $\mao{gcd}(w_j)=1$ 
	let $\D(w)=\CH(\{\vec m\})$ denote the Newton polytope of a generic 
	quasi-homogeneous polynomial $\sum_m c_{m}\,x_j^{m_j}$ of degree
	$D=\sum_jw_j$ regarded as a lattice polytope on 
	the intersection of the lattice~$\IZ^{d+1}$ with the hyperplane 
	$\vec m\cdot\vec w=D$. Since $w_j>0$ and $m_j\ge0$ the Newton 
	polytope $\D(w)$ has at most one interior point, namely $m_j=1$. 
	We call $w$ an \IP-weight-vector if $\D(w)$ is an \IP-polytope.
\end{definition}
\begin{remark}
	It follows already from a more general result of Lagarias and Ziegler 
	\cite{LZ91} that the number of \IP\ polytopes, and hence the number 
	of \IP\ weight vectors, is finite in fixed dimension. 
	A constructive proof has been found in terms of an efficient algorithm
	for the enumeration of all \IP-weight-vectors \cite{S,pwf}. 
\end{remark}
\begin{remark}
	Since all \IP-polygons are reflexive, linear relations of vertices 
	of \IP-simplices in two dimensions are always \IP-weight-vectors. 
	For $d\ge3$, on the other hand, the latter 
	form a proper subset of the set of all linear relations among 
	simplicial \IP-polytopes. The lowest degree examples of 
	non-\IP-weights of \IP\ simplices in 3, 4 and 5 dimensions
	are $\mathbb P_{1,5,6,8}[20]$, $\mathbb P_{1,1,5,5,8}[20]$
	and $\mathbb P_{2,2,2,3,3,11}[23]$, respectively.
	Table\,\ref{tab:wv3d} contains the complete list for $d=3$ in 
	boldface. 
	Note that $w_j\le D/2$ for all linear relations of \IP-simplicies 
	because for $w_0>w_1+\ldots+w_d$ 
	the simplex $\big\langle v_j\big\rangle$
	with $\sum w_jv_j=0$ has at least two interior lattice points, 
	namely 0 and $-v_0=\sum c_jv_j$ with 
	$c_0=\frac{2w_0}D-1$ and $c_j=\frac{2w_j}Dv_j$ for $j>0$, so that
	$\sum c_j=1$ and $c_j>0$;
	if $w_0=D/2$ then $-v_0$ is in the interior of a facet.
\del
cws.x -w3 4 99 -I
20  1 5 6 8 ipc=0 nip=1
21  1 4 7 9 ipc=0 nip=1
24  2 5 8 9 ipc=0 nip=1
28  1 5 8 14 ipc=0 nip=1
28  3 7 8 10 ipc=0 nip=1
30  4 7 9 10 ipc=0 nip=1
33  5 8 9 11 ipc=0 nip=1
36  3 7 8 18 ipc=0 nip=1
44  5 8 9 22 ipc=0 nip=1
#primepartitions=119430 #IPCpolys=104

20  1 1 5 5 8
23  2 2 3 5 11 ipc=0 nip=1
23  2 3 3 4 11 ipc=0 nip=1
25  1 3 3 8 10 ipc=0 nip=1
25  1 1 6 7 10 ipc=0 nip=1
26  1 1 5 8 11 ipc=0 nip=1
27  2 2 3 7 13 ipc=0 nip=1
27  2 3 4 5 13 ipc=0 nip=1
28  1 1 5 8 13 ipc=0 nip=1
28  1 4 4 6 13 ipc=0 nip=1
28  3 3 4 5 13 ipc=0 nip=1
28  1 1 5 9 12 ipc=0 nip=1
28  1 1 7 8 11 ipc=0 nip=1
28  1 2 7 7 11 ipc=0 nip=1
28  1 5 5 8 9 ipc=0 nip=1
28  2 3 7 8 8 ipc=0 nip=1
28  1 5 6 8 8 ipc=0 nip=1
28  3 3 7 7 8 ipc=0 nip=1
29  2 3 5 5 14 ipc=0 nip=1
29  3 3 4 5 14 ipc=0 nip=1
30  3 4 4 5 14 ipc=0 nip=1
30  1 1 7 9 12 ipc=0 nip=1
30  1 2 7 8 12 ipc=0 nip=1
30  2 2 5 10 11 ipc=0 nip=1
30  2 5 5 7 11 ipc=0 nip=1
30  1 4 7 9 9 ipc=0 nip=1

23  2 2 2 3 3 11
25  1 1 1 6 6 10
27  2 2 2 3 5 13
27  2 2 3 3 4 13
28  1 1 1 7 7 11
28  1 1 5 5 8 8
29  2 2 3 3 5 14
29  2 3 3 3 4 14

\enddel
\end{remark}

\del
\begin{table}\footnotesize
\begin{tabular}{||c|c||}\hline\hline
4&1 1 1 1\\
5&1 1 1 2\\
6&1 1 2 2\\
6&1 1 1 3\\
7&1 1 2 3\\
8&1 2 2 3\\
8&1 1 2 4\\
9&1 2 3 3\\
9&1 1 3 4\\
10&1 2 3 4\\
10&1 2 2 5\\
10&1 1 3 5\\
11&1 2 3 5\\
12&2 3 3 4\\
12&1 3 4 4\\
12&2 2 3 5\\
12&1 2 4 5\\
12&1 2 3 6\\
12&1 1 4 6\\
13&1 3 4 5\\
14&2 3 4 5\\
14&2 2 3 7\\
14&1 2 4 7\\
15&3 3 4 5\\
15&2 3 5 5\\
15&1 3 5 6\\
15&1 3 4 7\\
15&1 2 5 7\\
16&1 4 5 6\\
16&2 3 4 7\\
16&1 3 4 8\\
16&1 2 5 8\\
17&2 3 5 7\\
18&3 4 5 6\\
18&1 4 6 7\\
\hline\hline\end{tabular}\begin{tabular}{||c|c||}\hline\hline
18&2 3 5 8\\
18&2 3 4 9\\
18&1 3 5 9\\
18&1 2 6 9\\
19&3 4 5 7\\
20&2 5 6 7\\
20&3 4 5 8\\
\bf20&\bf  1 5 6 8\\
20&2 4 5 9\\
20&2 3 5 10\\
20&1 4 5 10\\
21&3 5 6 7\\
21&1 5 7 8\\
21&2 3 7 9\\
\bf21&\bf  1 4 7 9\\
21&1 3 7 10\\
22&2 4 5 11\\
22&1 4 6 11\\
22&1 3 7 11\\
24&3 6 7 8\\
24&4 5 6 9\\
\bf24&\bf  2 5 8 9\\
24&1 6 8 9\\
24&3 4 7 10\\
24&2 3 8 11\\
24&3 4 5 12\\
24&2 3 7 12\\
24&1 3 8 12\\
25&4 5 7 9\\
26&2 5 6 13\\
26&1 5 7 13\\
26&2 3 8 13\\
27&5 6 7 9\\
27&2 5 9 11\\
\bf28&\bf  3 7 8 10\\
\hline\hline\end{tabular}\begin{tabular}{||c|c||}\hline\hline
28&4 6 7 11\\
28&3 4 7 14\\
\bf28&\bf  1 5 8 14\\
28&1 4 9 14\\
\bf30&\bf  4 7 9 10\\
30&5 6 8 11\\
30&3 4 10 13\\
30&4 5 6 15\\
30&2 6 7 15\\
30&1 6 8 15\\
30&2 3 10 15\\
30&1 4 10 15\\
32&4 5 7 16\\
32&2 5 9 16\\
\bf33&\bf  5 8 9 11\\
33&3 5 11 14\\
34&4 6 7 17\\
34&3 4 10 17\\
36&7 8 9 12\\
\bf36&\bf  3 7 8 18\\
36&3 4 11 18\\
36&1 5 12 18\\
38&5 6 8 19\\
38&3 5 11 19\\
40&5 7 8 20\\
42&3 4 14 21\\
42&2 5 14 21\\
42&1 6 14 21\\
\bf44&\bf  5 8 9 22\\
44&4 5 13 22\\
48&3 5 16 24\\
50&7 8 10 25\\
54&4 5 18 27\\
66&5 6 22 33\\  &\\\hline\hline
\end{tabular}
\\
The lowest degree examples of non-\IP-confined
	\IP\ simplices in 3, 4 and 5 dimensions
	are $\mathbb P_{1,5,6,8}[20]$, $\mathbb P_{1,1,5,5,8}[20]$
	and $\mathbb P_{2,2,2,3,3,11}[23]$, respectively.
\caption{Linear relations of \IP\ simplices that are not \IP\ weight vectors.
	\label{tab:wv3d}\HS-259 }
\end{table}
\enddel

\begin{table*}{\footnotesize	
\def\TabHead{\hline\hline\end{tabular}\begin{tabular}{c|c||}\hline\hline
		$D$ & $\!\!w_1\ldots w_4\!\!$ \\\hline}
\begin{tabular}{||c|c||}\hline\hline
		$D$ & $\!\!w_1\ldots w_4\!\!$ \\\hline
4&1 1 1 1\\
5&1 1 1 2\\
6&1 1 2 2\\
6&1 1 1 3\\
7&1 1 2 3\\
8&1 2 2 3\\
8&1 1 2 4\\
9&1 2 3 3\\
9&1 1 3 4\\
10&1 2 3 4\\
10&1 2 2 5\\
10&1 1 3 5\\
11&1 2 3 5\\
12&2 3 3 4\\
12&1 3 4 4\\
\TabHead
12&2 2 3 5\\
12&1 2 4 5\\
12&1 2 3 6\\
12&1 1 4 6\\
13&1 3 4 5\\
14&2 3 4 5\\
14&2 2 3 7\\
14&1 2 4 7\\
15&3 3 4 5\\
15&2 3 5 5\\
15&1 3 5 6\\
15&1 3 4 7\\
15&1 2 5 7\\
16&1 4 5 6\\
16&2 3 4 7\\
\TabHead
16&1 3 4 8\\
16&1 2 5 8\\
17&2 3 5 7\\
18&3 4 5 6\\
18&1 4 6 7\\
18&2 3 5 8\\
18&2 3 4 9\\
18&1 3 5 9\\
18&1 2 6 9\\
19&3 4 5 7\\
20&2 5 6 7\\
20&3 4 5 8\\
\bf20&\bf  1 5 6 8\\
20&2 4 5 9\\
20&2 3 5 10\\
\TabHead
20&1 4 5 10\\
21&3 5 6 7\\
21&1 5 7 8\\
21&2 3 7 9\\
\bf21&\bf  1 4 7 9\\
21&1 3 7 10\\
22&2 4 5 11\\
22&1 4 6 11\\
22&1 3 7 11\\
24&3 6 7 8\\
24&4 5 6 9\\
\bf24&\bf  2 5 8 9\\
24&1 6 8 9\\
24&3 4 7 10\\
24&2 3 8 11\\
\TabHead
24&3 4 5 12\\
24&2 3 7 12\\
24&1 3 8 12\\
25&4 5 7 9\\
26&2 5 6 13\\
26&1 5 7 13\\
26&2 3 8 13\\
27&5 6 7 9\\
27&2 5 9 11\\
\bf28&\bf  3 7 8 10\\
28&4 6 7 11\\
28&3 4 7 14\\
\bf28&\bf  1 5 8 14\\
28&1 4 9 14\\
\bf30&\bf  4 7 9 10\\
\TabHead
30&5 6 8 11\\
30&3 4 10 13\\
30&4 5 6 15\\
30&2 6 7 15\\
30&1 6 8 15\\
30&2 3 10 15\\
30&1 4 10 15\\
32&4 5 7 16\\
32&2 5 9 16\\
\bf33&\bf  5 8 9 11\\
33&3 5 11 14\\
34&4 6 7 17\\
34&3 4 10 17\\
36&7 8 9 12\\
\bf36&\bf  3 7 8 18\\
\TabHead
36&3 4 11 18\\
36&1 5 12 18\\
38&5 6 8 19\\
38&3 5 11 19\\
40&5 7 8 20\\
42&3 4 14 21\\
42&2 5 14 21\\
42&1 6 14 21\\
\bf44&\bf  5 8 9 22\\
44&4 5 13 22\\
48&3 5 16 24\\
50&7 8 10 25\\
54&4 5 18 27\\
66&5 6 22 33\\  &\\\hline\hline
\end{tabular}}
\\[9pt]
\caption{\label{tab:wv3d}The 104 \IP-simplices in $d=3$ correspond to 
	95 transversal
	and 9 (boldface) non-\IP-weights.\hspace*{-248pt}
}
\vspace{-9pt}
\end{table*}

\del	The lowest degree examples of non-\IP-confined
	\IP\ simplices in 3, 4 and 5 dimensions
	are $\mathbb P_{1,5,6,8}[20]$, $\mathbb P_{1,1,5,5,8}[20]$
	and $\mathbb P_{2,2,2,3,3,11}[23]$, respectively.
\enddel

\begin{theorem}[Skarke \cite S] 
	The Newton polytope $\D(w)$ of an \IP\ weight vector 
	$\vec w=(w_0,\ldots,w_d)$ is always reflexive in $d\le4$ dimensions.
	More generally, for weight-matrices $W_{ij}$ with 
	$1\le i\le I$, $1\le j\le J$ and $D_i=\sum_jW_{ij}$ the Newton
	polytopes $\D(W)$ of generic polynomials of multi-degrees $D_i$
	are always reflexive whenenver $W$ is an \IP-weight-matrix 
	(i.e. when $\D(W)$ is an \IP\ polytope) if the dimension $d=J-I$
	of $\D(W)$ obeys  $d\le4$.
\\[3pt]	{\rm The first counterexample is $\vec w=(1,1,1,1,1,2)$ with 
	degree $D=7$, whose 5-dimensional 
	Newton polytope $\D(w)$ has one facet with lattice-distance 2
	from the interior point and hence is not reflexive.
	}
\end{theorem}

\begin{remark}The number of \IP\ weight vectors is $N_{I\!P}(1)=1$, 
	$N_{I\!P}(2)=3$, $N_{I\!P}(3)=95$ and $N_{I\!P}(4)= 184\,026$ 
	for $d\le4$. By extrapolation we expect that an enumeration might be
	possible in $d=5$  but certainly not above because the numbers become
	too large. For enumeration problems discussed in the next section it
	would hence be important to find an efficient generalization of
	Skarke's algorithm \cite{S} for \IP\ weight vectors with a bounded
	number of lattice points in their defining minimal simplex.
\end{remark}
\begin{definition}
	A weight vector $w$ is called transversal
	if a generic quasi-homogeneous polynomial
	$f_{\D(w)}=\sum_{m\in\D}\,c_m\prod_j x_j^{m_j}$ of degree $D$
	is transversal $df(\vec x)=0~\Rightarrow~\vec x=0$,
	i.e. the gradient only vanishes at the origin $x_j=0$.
	Due to the Euler formula $f(\vec x)=\sum q_j x_j\,\partial f/
	\partial x_j$ with $q_j=w_j/D$ the equation $f=0$ thus defines an
	isolated singularity \cite{Kreuzer:1992bi}.
\end{definition}

\begin{remark}
	Transversal weight vectors are a subclass of the \IP\ weight 
	vectors which concide with the reflexive weights for $d\le3$ and 
	imply reflexivity for $d=4$ \cite S. They have been enumerated
	for $d\le5$ \cite{nms,t4d,4fold,t5d} with
	$N_T(1)=1$, $N_T(2)=3$, $N_T(3)=95$, $N_T(4)=7555$ and
	$N_T(5)=1\,100\,055$. The reflexive weights are much more
	numerous than the transversal ones already in 4 dimension.
	Nevertheless,
	only 252\,933 of the transversal weights in 5 dimensions, about 23\%,
	are reflexive. The above example $(1,1,1,1,1,2)$ of a non-reflexive
	\IP\ weight is also the simplest example of a non-reflexive
	transversal weight.
\end{remark}

\begin{table}[h]	
\medskip
\begin{tabular}{||c||c|c|c|c||}
\hline\hline
        & rd=1 & rd=2 & rd=3 & rd=4\\ \hline\hline
d=1     & 1     & 3     & 7     & 54    \\ \hline       
d=2     &       & 16    & 328   & 230109        \\ \hline  
d=3     &       &       &4\,319 & 45\,986\,238  \\ \hline  
d=4     &&&& 473\,800\,776\\\hline
\hline	\end{tabular}  
\bigskip
\caption{Numbers of lattice polytopes with reflexive dimension $\le4$.
	\label{tab:RD}\HS-259 }\vspace{-9mm}
\end{table}

The complete algorithm for the enumeration of reflexive polytopes is
compiled in \cite{pwf}. A somewhat subtle detail is 
the possible choice of a sublattice on the side of the maximal polytope, 
which is constrained by the fact that the surviving lattice points need 
to span an \IP\ polytope on the sublattice. An implementation of the 
algorithm and of additional applications to combinatorics and algebraic 
geometry has been published in the program package PALP \cite{palp}. 
With some further extensions, which are implemeted in the current version 1.1,
this package has been used for all computations of results that are 
report in this note. In particular, we enumerated all lattice polytopes 
with reflexive dimension%
\footnote{~
	The program and many data are available at \cite{CYdata}.
	Basic descriptions of e{\bf x}tended/e{\bf x}perimental features
	are available
	with the help option {\tt -x}.
	In particular, facets are enumerated with {\,\tt poly.x\,\,-F} and
	{\,\tt class.x\,\,-A\,} creates lists of 
	affine normal froms.
	The third column of \tab tab:RD. can be obtained with

	{\scriptsize{\alltt \$ 
echo "4 1 1 1 1 /Z2: 0 1 0 1" > 3d.in; cws.x -w3 >> 3d.in
	}\vspace{-3pt}{\alltt \$ 
cws.x -c3 >> 3d.in; class.x -po 3d.ref 3d.in
	}\vspace{-3pt}{\alltt \$
class.x -b -pi 3d.ref | poly.x -fF |\,class.x -f 
					\hbox to 0pt{-A -po 3rd.2d\hss}
	}\vspace{-3pt}{\alltt	\$
class.x -B -pi 3rd.2d | poly.x -fF |\,class.x -f -A -po 3rd.1d
	}}
\\[-9pt]
	after fetching and compiling PALP, e.g. as 
	in (\ref{FetchCompile}) below.
}
$\mao{rd}(\D)\le4$, whose numbers are listed in \tab tab:RD..

\bigskip

\Section{
	\IP-confined polytopes}

The algorithm for the enumeration of reflexive polytopes is simplified
in $d\le4$ 
by the fact that Newton polytopes 
$\D(W)$ of \IP\ weight matrices are automatically
reflexive so that all maximal reflexive polytopes are of this form.
Since this is no longer true for $d\ge5$ we are lead to consider a
slightly larger class of \IP\ polytopes.
\begin{definition}				\label{IPCdual}
	Let $\widetilde \D:=\CH(\D^*\cap\IZ^d)$ be the convex hull of 
	the lattice points in~$\D^*$. An \IP\ polytope is called \IP-confined 
	(IPC) if $\widetilde \D$ is a lattice polytope. The polytope 
	$\widetilde{\widetilde\D}\subset\D$
	is called the IPC-closure of $\D$ and $\D$ is called IPC-closed if 
	$\D=\widetilde{\widetilde\D}$.
\end{definition}
\begin{remark}
	The operation $\D\to\widetilde\D$ extends the polar duality of
	reflexive polytopes to an involution on the set of IPC-closed
	polytopes. (Reflexive polytopes are obviously IPC-closed; according
	to Skarke's theorem \cite S the reverse is true for $d\le4$.)
\end{remark}
\begin{corollary}
	Each IP-confined polytope (and hence also each reflexive polytope)
	is a subpolytope of a maximal IPC-closed polytope $\D_M$. These
	maximal polytopes are of the form $\D_M=\CH(\D(W)\cap\G)$ where 
	$W$ is an \IP-weight-matrix with $\mao{rank}(W)\le d$ and 
	$\G\subset \IZ^d$ 
	is a sublattice of the natural lattice of the Newton polytope for
	which $\D(W)\cap\G$ still has an interior point.
\end{corollary}
%

\begin{proof}
	The proof is analogous to the one for the reflexive case \cite{pwf}
	and uses
	the generalized duality $\D\mapsto\widetilde\D$ of definition 
	\ref{IPCdual}. The polytopes $\widetilde\D$ that are 
	dual to the maximal IPC-closed polytopes are minimal in the set of
	IPC-closed polytopes, but not necessarily minimal as IP polytopes.
	For every minimal IP polytope $\D'\subseteq\widetilde\D_M$, however,
	$\widetilde\D_M$ is the IPC-closure, because otherwise $\D_M$ would
	not be maximal. A weight matrix $W$ with $\D_M=\CH(\D(W)\cap\G)$
	can hence be constructed by regarding $\D'$ as the convex hull
	of lower-dimensional IP-simplices $\S_i$ (which are IP-confined since
	the duals to the restrictions to the linear subspaces of $\S_i$
	are projections of $\D_M$, which provide confining IP polytopes 
	for the simplices $\S_i$). Each simplex in the decomposition can
	be chosen such that it contains at least two vertices of $\D'$
	that are not contained in a different simplex of the decomposition. 
	The maximal rank $d$ hence occurs if the maximal polytope is a 
	hypercube. 
\end{proof}

The numbers of reflexive polytopes in $d$ dimensions are the
diagonal entries of \tab tab:RD.. A double-exponential ansatz
$N_d\approx 2^{2^{d+1}-4}$ due to Skarke provides a good
fit to these data and would predict $N_5\approx1.2\cdot10^{18}$ and
$N_6\approx2.1\cdot10^{37}$ in 5 and 6 dimensions, respectively.
If this gives the correct order of magnitude an enumeration would already
be hopeless in 5 dimensions because the result could not be stored on
any existing medium (for the 4d-case we needed, on average, about 
20 byte, or five 32-bit integers, to store the data of a reflexive pair).
Bounds on the numbers of lattice points have been obtained in \cite{Nil07}.

While the case of 4-dimensional reflexive polytopes is of direct relevance
to algebraic geometry and string theory because the corresponding toric
hypersurfaces are Calabi-Yau 3-folds \cite{B}, higher-dimensional cases
are important because F-theory applications \cite{Denef} in physics
require (elliptically fibered) Calabi-Yau 4-folds, and, more
generally, reflexive polytopes in $n+r$ dimensions are related to complete 
intersection $n$-folds of codimension $r$ \cite{BBnef,BBsth,BN}. From
a pragmatical point of view, in both contexts the polytopes with a small 
number of points are the ones that we are most interested in because they 
encode the combinatorial data for manifolds with (moderately) small 
Picard numbers. This suggests a different setting for enumeration attempts
where we ask for the lists of reflexive polytopes with a certain (maximal)
number of points rather than a fixed dimension. Instead of descending
from maximal objects, which is not feasible for $d>4$, one might hence
add points to a list of minimal objects.
\begin{question}	\label{Qalgo}
Is there an efficient algorithm for the enumeration of weight vectors $\vec w$
that define IP-confined simplices $\sum_j w_jv_j=0$ with a bounded number of lattice points?
\end{question}
While we expect the number of IP-weight-vectors to be too large for $d>5$
for a complete enumeration, we also expect that the enumeration of all 
reflexive polytopes with less than, say, 12 or 15 lattice points should be
feasible if the questions \ref{Qalgo} has a positive answer. 
Here the natural objects of interest are IP-confined lattice
polytopes because every IPC-closed polytope, and hence every reflexive 
polytope, can be constructed by successively adding lattice points to the
minimal objects in a finite number of ways   as 
determined by the 
relevant weight vectors.

\bigskip

\Section{Statistics and numerical experiments}

\begin{figure*}[th]
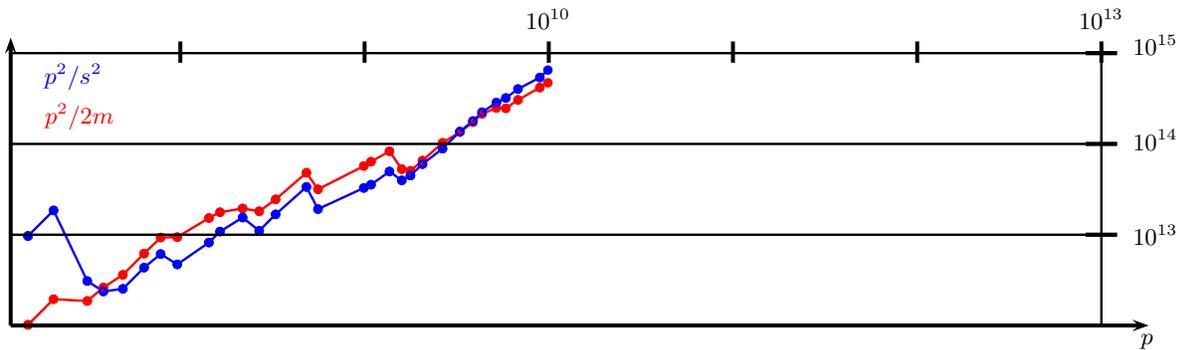

\TMstatistic\vspace{-2pt}
\caption{Predictions for 
	reflexive transversal Newton
	polytopes with $\le36$ $\ldots$ $\le65$ points.\label{fig:TM5}\HS-239 }
\vspace{2pt}
\end{figure*}

We now turn to statistical considerations that may provide additional
information in 
cases where a complete enumeration is not possible.
The first question we want to address is whether the property of self-duality
of reflexive lattice polytopes is random in the sense that the number of
selfdual objects can be estimated by a probabilistic calculation.


Let $N$ be the number of reflexive polytopes in fixed dimension $d$.
Then the number of (formal) duality assignments with $S$ self-mirrors, i.e. 
the number of involutions with $S$ fixed points, is $N\choose S$
for the choice of $S$ self-mirrors times 
$(N-S-1)\cdot(N-S-3)\cdot\ldots \cdot3\cdot 1$ for the possible selections of
the remaining dual pairs, hence
%
%
\BE
	n_S=
	\frac{N!}{{S!}\cdot2^{\frac{N-S}2}\,(\frac{N-S}2)!}\qquad
	\hbox{with}\quad S-N\in2\IZ.
\EE
\del
For large $N$ we can use
Stirling's formula 
\BE	\textstyle
	N!=\sqrt{2\p N} N^Ne^{-N+\frac\th{12N}}
	=\sqrt{2\p N}\left(\frac Ne\right)^n\left(1+\mathcal O(\frac1N)\right)
\EE
with $0<\th<1$
\BE	\textstyle
	\frac{\6\log(n_{2s})}{\6N}\approx \frac1{\sqrt{2N}}
	\left(1+\frac sN+\mathcal O(\frac{s^2}{N^2})\right)
\EE
hence with $Z(N)=\sum_{0\le s\le N/2}\,2s\,n_{2s}$
\BE	\textstyle
	\6_N\log Z\approx\frac1{\sqrt{2N}}
			\left(1+\frac{\langle S\rangle}N\right)
\EE
\BE	\textstyle
	\6_S n_{S}=
\EE
\BE
	n_S=\frac {(\frac N2)!\;2^{S/2}} {S!\;(\frac{N-S}2)!}	WRONG
\EE
where we assumed that $N$, and hence also the number $S$
of self-mirror polytopes, is even. 	
This yields the expectation value  
\BE	\textstyle
	\langle S\rangle={\sum 2s\,n_{2s}\over\sum n_{2s}}\approx
	{\sum 2s\,(2/N)^s/(2s)!\over\sum(2/N)^s/(2s)!}
	=-2N\6_N\log(\textstyle{\sum{(2/N)^s\over(2s)!}})
	\approx-2N\6_N\log\cosh\8{\sqrt{2\0N}}
\EE								\enddel
Let $Z_N=\sum_{S\le N} n_S$ be the total number of involutions.
The asymptotic expansion for large $N$ 
can be derived from the generating function \cite[section 3.8]{genfun}
\BE
	e^{x+\frac12x^2}=\sum_{N\ge0}\frac{Z_N}{N!}x^N.
\EE
Since it happens $Z_{N-1}$ times for the $Z_N$ involutive 
permutations of $N$ objects that a given object is a fixed point,
we obtain the following formula for the expectation value \cite{ILugo}
\BE	\textstyle
	\big\langle S\big\rangle=\frac1{Z_N}\sum_Sn_S\,S=N\frac{Z_{N-1}}{Z_N}.
\EE
Assuming a uniform probability distribution of involutions 
we thus expect
\BE	\textstyle
	\big\langle S\big\rangle=\sqrt N-\frac12+\frac13\frac1{\sqrt N}
			+{\mathcal O}(\frac1N)
\EE
selfdual polytopes for 
large $N$.
As shown in \tab tab:SD. this roughly explains the size but 
underestimates the correct numbers for $d\le4$.

\begin{table}[h]
\BC\begin{tabular}{||c||c|c|c|c||}\hline\hline
$d$	& 1&2&3&4\\\hline\hline
$N_R(d)$& 1&16&4319&473800776\\\hline
$N_{\rm selfdual}$&1&4&79&41710	\\\hline
$\big\langle S\big\rangle$ 
	& 1 & 3.6 & 65.2 & 21766.5	\\\hline\hline
\end{tabular}\\[7pt]
\caption{
	Numbers of selfdual polytopes and 
	probabilistic 
	expectations
	.\label{tab:SD}\HS-259
	}	\vspace{-18pt}
\EC\end{table}



\del
\BE	\hspace*{-9pt}
	\frac{\sum_{s\ge0} \,2s\,\, n_{2s}}{\sum_{s\ge0}n_{2s}}
	\approx 2N\6_N\log\sum_{s\ge0} {N^s\over(2s)!}=2N\6_N\log\cosh\sqrt N
	=\sqrt N\tanh\sqrt N\approx\sqrt N
\EE
With $s=S/2\in\IZ_{\ge0}$ this yields the expectation value  
\BE	\langle S\rangle={\sum 2s\,n_{2s}\over\sum n_{2s}}\approx
	{\sum 2s\,(2/N)^s/(2s)!\over\sum(2/N)^s/(2s)!}
	=-2N\6_N\log(\textstyle{\sum{(2/N)^s\over(2s)!}})
	\approx-2N\6_N\log\cosh\8{\sqrt{2\0N}}
\EE
For large $N$ we thus find
\BE	\langle S\rangle\approx 2N\tanh\8{\sqrt{2\0N}}\;\8{\2\sqrt{N\02}
		{1\0N^2}}=1/N
\EE
							\enddel

\del
The numbers $N_2=16$, $N_3=4\,319$ and $N_4=473\,800\,776$ of reflexive 
polytopes in 2, 3 and 4 dimensions show that a complete enumeration of the 
5-dimensional case is hopeless at present already because no existing 
harddisk would be large enough to store the data. One may try to estimate 
their number by extrapolation. The double-exponential ansatz
$N_d\approx 2^{2^{d+1}-4}$ due to Harald Skarke provides a good 
fit to the known data and would predict $N_5\approx1.2\cdot10^{18}$ and 
$N_6\approx2.1\cdot10^{37}$ in 5 and 6 dimensions, respectively. 
For a random choice we would expect the first mirror pairs to show up 
at $\sqrt{N_d}$, i.e. at $10^9$ polytopes ... this is about twice the upper
limit of the maximal content of one PALP database with
present computer architecture and amounts to 30GB of data. 
Appears to be hard to get a reasonable statistics.
\enddel

As this result seems to indicate that a statistical approach makes sense we
now want to use similar considerations for predicting the number of
reflexive polytopes in 5 dimensions on the basis of incomplete lists.
In a random set of $p>\sqrt N$ reflexive polytopes in fixed dimension
the formula $\langle S\rangle\approx
\sqrt N$ implies that we expect  
$s\approx\sqrt N \cdot p/N$ self-mirror polytopes.
This leads to the prediction $N\approx (p/s)^2$ if we find $s$ selfdual
polytopes in the sample.
Similarly, if we ignore the relatively small number of self-mirrors 
for large $N$ and count the number $m$ of mirrors pairs in a sample of
$p$ polytopes we obtain the prediction $N\approx p^2/(2m)$. 
If we increase the size of a random sample we expect $s$ to grow linearly 
and $m$ quadratically with $p$, and more precisely $s\approx\sqrt{2m}$.

Actually, the formula $N\approx p^2/(2m)$
has been used already serval years ago when we enumerated the reflexive
polytopes in 4 dimensions \cite{c4d}, which required two years of program
improvements and computation time after the 3-dimensional case \cite{c3d}.
We thus could check the sufficiency of the implemented data structures for 
the storage of the result at an early stage of the project.
The starting point of the calculation was the list of 308 weight matrices
(206 weight vectors and 102 matrices with $2\le\mao{rank}\le4$) of
maximal reflexive polytopes $\D_M$ for which the reflexive subpolytopes 
were computed in the order of an increasing number of lattice points.
The first Newton polytope in this list is defined by the weight vector 
$(3,3,4,4,10)$ with degree 24. It has 47 lattice points and 6 vertices.
After fetching and compiling PALP its subpolytopes can be computed 
with the following commands,
\begin{equation}
\hspace*{-34mm}\hbox to 5cm {
\vbox{		\tiny\label{FetchCompile}\begin{verbatim}
$ Bdir=$HOME/bin               # directory for binary files (check $PATH)
$ Wdir=$PWD                    # working directory
$ cd /tmp 
$ wget hep.itp.tuwien.ac.at/~kreuzer/CY/palp/palp-1.1.tar.gz     # fetch
$ gunzip palp-*.tar.gz; tar -xvf palp-*.tar; cd palp;            # unpack
$ make; mv *.x $Bdir; cd $Wdir                         # and compile PALP
$
$ echo "24 3 3 4 4 10" | class.x -f -po /tmp/zbin.47
...
800kR-1658 11MB 4215kIP 1042kNF-23k 9_46 v17r17 f28r27 85b24 120s 120u 25n
24 3 3 4 4 10 R=798878 +1658sl hit=0 IP=4215623 NF=1042005 (23261)
... 1181m+14s 9851469b  u7 pp/2m=2.68615e+08 pp/ss=3.25615e+09
\end{verbatim}\hss}}
\end{equation}
Within two minutes we thus obtain
the values 269 million for $p^2/2m$ with $p\!=\!798878$ and $m\!=\!1181$
and 3.25 billion for $p^2/s^2$ with $s\!=\!14$ on a standard 3GHz PC 
(cf.\;the last output line) and hence a good approximation of the correct value 473\,800\,776 with
a production rate of about $7000$ polytopes per second (back in 1998 the
CPU time was almost 1 hour).

					\del
\begin{verbatim}
$ echo "15 2 2 3 3 5" | class.x -f -po /tmp/zbin56	# M:56 9 N:8 6
...
11202kR-3930 167MB 77MIP 14MNF-82k 7_55 v22r21 f32r29 61b30 47m 47u 8n
15 2 2 3 3 5 R=11198394 +3930sl hit=0 IP=77167925 NF=14872004 (82172)
... 11198394+3930sl 34804m+71s 160699062b  u13 pp/2m=9.27455e+08 pp/ss=7.40524e+09

$ echo "4 1 1 1 1 0 0  6 0 0 1 1 2 2" |class.x -f -po /tmp/zbin # M:57 8 N:9 6
...
17791kR-2445 238MB 105MIP 27MNF-255k 7_55 v25r23 f36r36 46b31 68m 68u 15n
4 1 1 1 1 0 0 6 0 0 1 1 2 2 R=17789212 +2445sl hit=0 IP=105672169 NF=27019369 (255903)
... 17789212+2445sl 795011m+1482s 235323323b  u23 pp/2m=1.86804e+07 pp/ss=4.55277e+06

echo "14 2 2 3 3 4" | class.x -f -po /tmp/zbin		# M:57 10 N:10 7

wget http://quark.itp.tuwien.ac.at/~kreuzer/W/4dTransWH.gz;
gunzip 4dTransWH.gz ; mv 4dTransWH 4dT     # cws.x -w4 -t > 4dT
cws.x -w3 -t > 3dT ; cws.x -w2 -t > 2dT ; cws.x -w1 -t > 1dT

wget http://quark.itp.tuwien.ac.at/~kreuzer/W/5dTransWH.all.gz

for((i=10; i<=27; i++));
     do echo "$i:"; cat tc52* | grep "M:$i " >> tcm10-27.in; done

echo "24 3 3 4 4 10" | class.x -f -po /tmp/zbin
Writing zbin: 798878+1658sl 1181m+14s 9851469b  u7 pp/2m=2.68615e+08 pp/ss=3.25615e+09

echo "24 3 3 4 4 10" | class.x -f -po zbin
798878+1658sl 1181m+14s 9851469b
pp/2m=2.68615e+08
pp/ss=3.25615e+09

for((i=28;i<=30;i++));do echo "$i:"; cat tc52*|grep "M:$i ">>tcm$i.in; done

beauty/tmp> echo "24 3 3 4 4 10" | class.x -f -po /tmp/zbin
Fri Mar 21 15:40:56 2008
100kR-415 1MB 453kIP 111kNF-2k 7_47 v16r15 f25r22 85b21 12s 12u 2n
200kR-1168 2MB 857kIP 227kNF-5k 6_37 v16r16 f26r24 85b21 24s 24u 5n
800kR-1658 11MB 4215kIP 1042kNF-23k 9_46 v17r17 f28r27 85b24 121s 121u 25n
24 3 3 4 4 10 R=798878 +1658sl hit=0 IP=4215623 NF=1042005 (23261)
Writing /tmp/zbin: 798878+1658sl 1181m+14s 9851469b  u7 pp/2m=2.68615e+08 pp/ss=3.25615e+09 done: 1s
Fri Mar 21 15:42:58 2008

brane/palp> echo "24 3 3 4 4 10" | class.x -f -po /tmp/zbin
Fri Mar 21 15:42:01 2008
100kR-415 1MB 453kIP 111kNF-2k 7_47 v16r15 f25r22 85b21 15s 15u 3n
200kR-1168 2MB 857kIP 227kNF-5k 6_37 v16r16 f26r24 85b21 30s 30u 6n
300kR-1024 4MB 1327kIP 345kNF-7k 8_40 v16r16 f27r25 85b21 47s 47u 10n
700kR-1221 9MB 3553kIP 892kNF-17k 6_42 v17r17 f28r27 85b24 127s 126u 26n
800kR-1643 11MB 4208kIP 1041kNF-22k 7_46 v17r17 f28r27 85b24 150s 150u 31n
800kR-1658 11MB 4215kIP 1042kNF-23k 9_46 v17r17 f28r27 85b24 150s 150u 31n
24 3 3 4 4 10 R=798878 +1658sl hit=0 IP=4215623 NF=1042005 (23261)
Writing /tmp/zbin: 798878+1658sl 1181m+14s 9851469b  u7 pp/2m=2.68615e+08 pp/ss=3.25615e+09 done: 0s
Fri Mar 21 15:44:31 2008
\end{verbatim}
		\enddel

For 5 dimensions it is, of course, much harder to get a reliable statistics.
With an expectation of $N\approx 10^{18}$ according to Skarke's guess, the
storage of $\sqrt N$ polytopes would already required some 30GB of disk space
so that we can only go above that value by 1-2 orders of magnitude with 
currently available hardware. Nevertheless, it should be possible
to either verify that the number is not much smaller or to get a reasonable
prediction if it is.
For a first attempt we defined data samples in terms of the Newton polytopes
of transversal reflexive weight vectors%
\footnote{~ cf.
\tt http://quark.itp.tuwien.ac.at/{\tiny$^\sim$}kreuzer/W/5dTR*
.}
ordered according to increasing numbers of lattice points. For the increasing
series data samples consisting of all reflexive subpolytopes of transversal 
Newton polytopes with $\le36$ \ldots $\le65$ points 
the predictions are plotted in \fig fig:TM5.

The result is obviously inconclusive, but certainly compatible with 
the guess $10^{18}$. Unfortunately the bias of the a priori independent
predictions $p^2/2m$ and $p^2/s^2$ is strongly correlated in our data samples.
For the $9.025\cdot10^{9}$ polytopes of the largest sample
with $m=86323$ and $s=354$, which occupies 239 GB of disk space, we find
$p^2/2m\approx4.7 \cdot10^{14}$ and $p^2/s^2\approx6.5 \cdot10^{14}$.

\bigskip

\Section{Conclusions}

In this note we determined all lattice polytopes with reflexive dimension 
$\mao{rd}\le4$ and discussed enumeration problems and algorithmic aspects
with applications to algebraic geometry and string theory.
We pointed out the need for an efficient algorithm for the enumeration
of IP weight vectors $w$ with a bounded number of lattice points in 
the convex hull of 
the simplex 
defined by the linear relations $\sum w_jv_j=0$. Such an algorithm could 
be used for the enumeration of reflexive polytopes with fixed number of
points rather than fixed dimension.

We introduced the concept of \IP-confined polytopes, which are a subclass
of \IP\ polytopes, and extended the polar duality of reflexive polytopes
to IPC-closed polytopes. Maximal IPC-closed polytopes contain all reflexive
polytopes in arbitrary dimensions and hence lead to a simplification of
the classification program. In turn, we pointed out the existence of 
\IP-simplices that are not \IP-confined and enumerated them for the case
of 3 dimensions. A constructive classification of such
simplices in higher dimensions is another interesting open problem.

We suggested a statistical
approach to the enumeration of reflexive polytopes which should at least
allow us to obtain probabilistic lower bounds, depending essentially on
the size of the available hard-disks for storage of the data.
As a first attempt we constructed a data-base containing about $9\cdot10^9$
pairs of reflexive 5-dimensional polytopes, which can also be used to 
produce incomplete lists of polytopes with reflexive dimension 5.
Studies of correlations of polytope data like $f$-vectors (numbers of
faces) can thus be initiated and may be useful for selecting appropriate 
data samples for statistical applications. 

\del
In the remaining months before the workshop 
``Information-theoretic aspects of integer-point enumeration in polyhedra''
{\footnotesize\verb+http://www.bio-complexity.com/ITSL/ITSL_index.html+}
we plan to create a complementary
sample based on weight matrices of rank two and to attempt an enumeration
of all IP weight vectors for dimension 5, which would enable
the enumeration of all maximal/minimal pairs of IPC-complete polytopes, 
and presumably also of the maximal/minimal pairs of reflexive polytopes.
More detailed data will be made available on the internet \cite{inprog}.
\enddel

\vspace{-2pt}
\def\refname{{\Fsec References \hspace*{165pt}}}
\end{document}

--------------------------------------------------------
11  1 1 1 1 1 6	
5 8  Points of P
   -1    1    0    0    0    0    0    0
   -1    0    1    0    0    0    0    0
   -1    0    0    1    0    0    0    0
   -1    0    0    0    1    0    0    0
   -6    0    0    0    0    1    0   -1
6 5  Equations of P		eval on -e5
  10  -1  -1  -1  -1     1
  -1  -1  -1  10  -1     1
  -1  -1  10  -1  -1     1
  -1  10  -1  -1  -1     1
  -1  -1  -1  -1  -1     1
  -6  -6  -6  -6   5     6
--------------------------------------------------------
11  1 2 2 2 2 2		not IP
5 8  Points of P
   -2    1    0    0    0    0   -1    0
   -2    0    1    0    0    0   -1    0
   -2    0    0    1    0    0   -1    0
   -2    0    0    0    1    0   -1    0
   -2    0    0    0    0    1   -1    0
6 5  Equations of P		eval on (-1,-1,-1,-1,-1)
   9  -2  -2  -2  -2     2	:: 0<= 1 + 2   non-IP
  -2  -2  -2   9  -2     2	:: 0<= 1 + 2 
  -2  -2   9  -2  -2     2	:: 0<= 1 + 2 
  -2   9  -2  -2  -2     2	:: 0<= 1 + 2 
  -1  -1  -1  -1  -1     1	:: 0<= 5 + 1 
  -2  -2  -2  -2   9     2	:: 0<= 1 + 2 

II=6 ; typeset -i deg ; 
for((i=2;i<=II;i++)) ; do 
for((j=1;j<=i;j++)) ; do 
for((k=1;k<=j;k++)) ; do 
for((l=1;l<=k;l++)) ; do 
	deg=$i+$j+$k+$l ; echo "$deg $i $j $k $l" ;done;done;done;done|
	cws.x -S|grep -v type | grep -v N:| poly.x -f| sort -u

\del
\begin{code}~\\[-7mm]{\footnotesize\begin{verbatim}
Bdir=$HOME/bin                  # directory for binary files (check $PATH)
Wdir=$PWD                       # working directory
cd /tmp ; wget hep.itp.tuwien.ac.at/~kreuzer/CY/palp/palp-1.1.tar.gz
gunzip palp-*.tar.gz; tar -xvf palp-*.tar; cd palp;make;mv *.x $Bdir;cd $Wdir

wget http://quark.itp.tuwien.ac.at/~kreuzer/W/4dTransWH.gz;
gunzip 4dTransWH.gz ; mv 4dTransWH 4dT     # cws.x -w4 -t > 4dT
cws.x -w3 -t > 3dT ; cws.x -w2 -t > 2dT ; cws.x -w1 -t > 1dT

wget http://quark.itp.tuwien.ac.at/~kreuzer/W/5dTransWH.all.gz

for((i=10; i<=27; i++));
     do echo "$i:"; cat tc52* | grep "M:$i " >> tcm10-27.in; done

echo "24 3 3 4 4 10" | class.x -f -po /tmp/zbin
Writing zbin: 798878+1658sl 1181m+14s 9851469b  u7 pp/2m=2.68615e+08 pp/ss=3.25615e+09

echo "24 3 3 4 4 10" | class.x -f -po zbin
798878+1658sl 1181m+14s 9851469b
pp/2m=2.68615e+08
pp/ss=3.25615e+09

for((i=28;i<=30;i++));do echo "$i:"; cat tc52*|grep "M:$i ">>tcm$i.in; done

beauty/tmp> echo "24 3 3 4 4 10" | class.x -f -po /tmp/zbin
Fri Mar 21 15:40:56 2008
100kR-415 1MB 453kIP 111kNF-2k 7_47 v16r15 f25r22 85b21 12s 12u 2n
200kR-1168 2MB 857kIP 227kNF-5k 6_37 v16r16 f26r24 85b21 24s 24u 5n
800kR-1658 11MB 4215kIP 1042kNF-23k 9_46 v17r17 f28r27 85b24 121s 121u 25n
24 3 3 4 4 10 R=798878 +1658sl hit=0 IP=4215623 NF=1042005 (23261)
Writing /tmp/zbin: 798878+1658sl 1181m+14s 9851469b  u7 pp/2m=2.68615e+08 pp/ss=3.25615e+09 done: 1s
Fri Mar 21 15:42:58 2008

brane/palp> echo "24 3 3 4 4 10" | class.x -f -po /tmp/zbin
Fri Mar 21 15:42:01 2008
100kR-415 1MB 453kIP 111kNF-2k 7_47 v16r15 f25r22 85b21 15s 15u 3n
200kR-1168 2MB 857kIP 227kNF-5k 6_37 v16r16 f26r24 85b21 30s 30u 6n
300kR-1024 4MB 1327kIP 345kNF-7k 8_40 v16r16 f27r25 85b21 47s 47u 10n
700kR-1221 9MB 3553kIP 892kNF-17k 6_42 v17r17 f28r27 85b24 127s 126u 26n
800kR-1643 11MB 4208kIP 1041kNF-22k 7_46 v17r17 f28r27 85b24 150s 150u 31n
800kR-1658 11MB 4215kIP 1042kNF-23k 9_46 v17r17 f28r27 85b24 150s 150u 31n
24 3 3 4 4 10 R=798878 +1658sl hit=0 IP=4215623 NF=1042005 (23261)
Writing /tmp/zbin: 798878+1658sl 1181m+14s 9851469b  u7 pp/2m=2.68615e+08 pp/ss=3.25615e+09 done: 0s
Fri Mar 21 15:44:31 2008

\end{verbatim}}\end{code}


					\enddel

\begin{minipage}{\columnwidth}
\includegraphics[width=3.81cm]{Bild}
\captionof{figure}[Bild (für Abbildungsverzeichnis)]{Bild bla blub (ist etwas länger)}
\label{fig:Bild}
\end{minipage}